\documentclass[12pt]{amsart}
\usepackage{amssymb, amsmath}
\newtheorem{defi}{Definition}[section]
\newtheorem{lema}{Lemma}[section]
\newtheorem{teo}{Theorem}[section]
\newtheorem{cor}{Corollary}[section]
\newcommand{\mb}[1]{\mathbb{#1}}
\newcommand{\mc}[1]{\mathcal{#1}}
\begin{document}
\title[Infinite multiplicity]{Infinite multiplicity of abelian subalgebras in free group
subfactors}
\author{Marius B. \c Stefan}
\address{UCLA Mathematics Department, Los Angeles, CA 90095-1555}
\email{stefan@math.ucla.edu}
 \abstract \noindent We obtain an estimate of
Voiculescu's (modified) free entropy dimension for generators of a
$\mbox{\!I\!I}_1$-factor $\mc{M}$ with a subfactor $\mc{N}$
containing an abelian subalgebra $\mc{A}$ of finite multiplicity.
It implies in particular that the interpolated free group
subfactors of finite Jones index do not have abelian subalgebras
of finite multiplicity or Cartan subalgebras.
\endabstract
\maketitle
\section{Introduction}\label{s1}
Cartan subalgebras arise naturally in the classical group measure
space construction. Thus, if $\alpha$ is a free action of a
discrete countable group $\Gamma$ on a measure space $(X,\mu)$,
then the cross-product von Neumann algebra
$L^\infty(X,\mu)\times_\alpha\Gamma$ contains a copy of
$L^\infty(X,\mu)$ as a Cartan subalgebra. More generally, a Cartan
subalgebra of a von Neumann algebra $\mc{P}$ is a maximal abelian
$*$-subalgebra of $\mc{P}$ whose normalizer generates $\mc{P}$
(regular MASA) and which is the range of a normal conditional
expectation (\cite{17}, \cite{3}). D. Voiculescu defined
(\cite{14}, \cite{15}) an original concept of (modified) free
entropy dimension $\delta_0$ and proved (\cite{15}) that
$\delta_0$ of any finite system of generators of a von Neumann
algebra which has a regular diffuse hyperfinite $*$-subalgebra
(regular DHSA) is $\leq 1$. This answered in the negative the
longstanding open question of whether every separable
$\mbox{\!I\!I}_1$-factor contains a Cartan subalgebra since the
free group factors $\mc{L}(\mb{F}_n)$ (von Neumann algebras
generated by the left regular representations $\lambda
:\mb{F}_n\rightarrow\mc{B}(l^2(\mb{F}_n))$, $2\leq n\leq\infty$)
have systems of generators with $\delta_0>1$. D. Voiculescu's
result about the absence of Cartan subalgebras in free group
factors was extended by L. Ge (\cite{4}) and K. Dykema (\cite{2})
who showed that these factors do not have abelian subalgebras of
multiplicity one and  of finite multiplicity, respectively. We
mention that if $\mc{A}$ is a Cartan subalgebra in a
$\mbox{\!I\!I}_1$-factor $N$, then $(\mc{A}\cup J\mc{A}J)''$ is a
MASA in $\mc{B}(L^2(\mc{N},\tau ))$ (\cite{3}, \cite{9}), hence
$\mc{A}$ is in particular an abelian subalgebra of multiplicity
one.

The interpolated free group factors $\mc{L}(\mb{F}_t)$
$(1<t\leq\infty )$ were introduced independently by K. Dykema
(\cite{2}) and F. R\u adulescu (\cite{10}) as a continuation of
the discrete series $ \mc{L}(\mb{F}_n)$, $2\leq n\leq\infty$. We
prove (Corollary \ref{c2}) that the subfactors of finite Jones
index in the interpolated free group factors do not have abelian
subalgebras of finite multiplicity either. The result is a
consequence of the estimate of (modified) free entropy dimension
(Theorem \ref{t3}) $\delta_0(x_1,\ldots ,x_m)\leq 2r+2v+3$, where
$x_1,\ldots ,x_m$ are self-adjoint generators of the
$\mbox{\!I\!I}_1$-factor $\mc{M}$, $r$ is the integer part of the
Jones index of $\mc{N}$ in $\mc{M}$ and $v$ is the multiplicity of
an abelian subalgebra $\mc{A}$ in $\mc{N}$.

Schreier's Theorem describes all subgroups of finite index $k$ in
the free group $\mb{F}_n$: any such subgroup is isomorphic to the
free group $\mb{F}_{1+k(n-1)}$. A von Neumann algebra analogue of
the fact that $\mb{F}_{1+k(n-1)}$ can be embedded with finite
index $k$ in $\mb{F}_n$ was proved by F. R\u adulescu (\cite{10}):
$\mc{L}(\mb{F}_{1+\lambda^{-1}(t-1)})$ can be embedded in
$\mc{L}(\mb{F}_t)$ with finite index $\lambda^{-1}$
$\forall\,1<t\leq \infty$ $\forall\,\lambda ^{-1}\in \lbrace
4\cos^2\frac{\pi}{k}:k\geq 3\rbrace$. On the other hand, at the
von Neumann algebra level, with $\mc{L}(\mb{F}_n)$ instead of
$\mb{F}_n$, it is no longer known whether Schreier's Theorem is
still true. However, two properties are preserved when passing to
free group subfactors of finite index: Haagerup approximation
property (\cite{6}) and primeness (\cite{11}) i.e., the
indecomposability as tensor product of type
$\mbox{\!I\!I}_1$-factors. Our result about the absence of abelian
subalgebras of finite multiplicity (and thus, of Cartan
subalgebras) is a third property that seems to support the
Schreier conjecture for free group subfactors.

We recall some results from D. Voiculescu's free probability
theory (\cite{13}, \cite{14}, \cite{15}) for the reader's
convenience. If $\mc{M}$ is a $\mbox{\!I\!I}_1$-factor with its
unique faithful normalized trace $\tau$ then $||x||_s=\tau
((x^*x)^\frac{s}{2})^\frac{1}{s}$ ($1<s<\infty$) denotes the
$s$-norm of $x\in\mc{M}$, $L^2(\mc{M},\tau )$ denotes the
completion of $\mc{M}$ with respect to the $2$-norm, and
$\mc{M}\subset\mc{B}(L^2(\mc{M},\tau ))$ is the standard
representation of $\mc{M}$. For an integer $c\geq 1$ let
$\mc{M}_c(\mb{C})$ and $\mc{M}_c^{sa}(\mb{C})$ be the set of all
$c\times c$ complex matrices and respectively, of all $c\times c$
complex self-adjoint matrices. Let further $\mc{U}_c(\mb{C})$ be
the unitary group of $\mc{M}_c(\mb{C})$, $\tau_c$ be the unique
normalized trace on $\mc{M}_c(\mb{C})$, and $||\cdot
||_e=\sqrt{c}||\cdot ||_2$ be the euclidian norm on
$\mc{M}_c(\mb{C})$. The free entropy of $x_1,\ldots
,x_m\in\mc{M}^{sa}$ in the presence of $x_{m+1},\ldots
,x_{m+n}\in\mc{M}^{sa}$ is defined in terms of sets of matricial
microstates $\Gamma_R((x_i)_{1\leq i\leq m}:(x_{m+j})_{1\leq j\leq
n};p,c,\epsilon )\subset (\mc{M}_c^{sa}(\mb{C}))^m$. The set
$\Gamma_R$ of matricial microstates corresponding to integers
$c,p\geq 1$ and to $\epsilon >0$ consists in $m$-tuples
$(A_i)_{1\leq i\leq m}$ of $c\times c$ self-adjoint matrices such
that there exists an $n$-tuple $(A_{m+j})_{1\leq j\leq
n}\in(\mc{M}_c^{sa}(\mb{C}))^n$ with the properties
\[|\tau (x_{i_1}\ldots
x_{i_l})-\tau_c (A_{i_1}\ldots A_{i_l})|<\epsilon\,,\,\,
||A_k||\leq R
\]
for all $1\leq i_1,\ldots ,i_l\leq m+n$, $1\leq l\leq p$, $1\leq k\leq m+n$.
One defines then successively:
\begin{eqnarray}&&\chi_R((x_i)_{1\leq i\leq m}:(x_{m+j})_{1\leq j\leq
n};p,c,\epsilon )\\\nonumber &&\hspace{.5 cm}
=\log\mbox{vol}_{mc^2}(\Gamma_R((x_i)_{1\leq i\leq
m}:(x_{m+j})_{1\leq j\leq n};p,c,\epsilon )),
\end{eqnarray}
\begin{eqnarray}&&\chi_R((x_i)_{1\leq i\leq m}:(x_{m+j})_{1\leq j\leq
n};p,\epsilon )\\\nonumber &&\hspace{.5 cm}
=\limsup_{c\rightarrow\infty}\left(\frac{1}{c^2}\chi_R((x_i)_{1\leq
i\leq m}:(x_{m+j})_{1\leq j\leq n};p,c,\epsilon )+\frac{m}{2}\log
c\right),
\end{eqnarray}
\begin{eqnarray}&&\chi_R((x_i)_{1\leq i\leq m}:(x_{m+j})_{1\leq j\leq
n})\\\nonumber &&\hspace{.5
cm}=\inf_{p,\epsilon}\chi_R((x_i)_{1\leq i\leq m}:
(x_{m+j})_{1\leq j\leq n};p,\epsilon ),
\end{eqnarray}
\begin{eqnarray}&&\chi((x_i)_{1\leq i\leq m}:(x_{m+j})_{1\leq j\leq
n})\\\nonumber &&\hspace{.5 cm}=\sup_R\chi_R((x_i)_{1\leq i\leq
m}:(x_{m+j})_{1\leq j\leq n})
\end{eqnarray}
(we denoted by $\mbox{vol}_{mc^2}(\cdot )$ the Lebesgue measure on
$(\mc{M}_c^{sa}(\mb{C}))^m\simeq\mb{R}^{mc^2}$). The resulting
quantity $\chi ((x_i)_{1\leq i\leq m}:(x_{m+j})_{1\leq j\leq n})$
is the free entropy of $(x_i)_{1\leq i\leq m}$ in the presence of
$(x_{m+j})_{1\leq j\leq n}$ or if $n=0$, the free entropy $\chi
(x_1,\ldots ,x_m)$ of $(x_i)_{1\leq i\leq m}$. The free entropy of
$(x_i)_{1\leq i\leq m}$ in the presence of $(x_{m+j})_{1\leq j\leq
n}$ is equal to the free entropy of $(x_i)_{1\leq i\leq m}$ if
$\{x_{m+1},\ldots ,x_{m+n}\}\subset\{x_1,\ldots ,x_m\}''$. Also,
the free entropy of a single self-adjoint element $x$ is (where
$\mu$ denotes the distribution of $x$): $$\chi
(x)=\frac{3}{4}+\frac{1}{2}\log 2\pi +\int\int\log |s-t|d\mu
(s)d\mu (t).$$ An element $x\in \mc{M}$ is a semicircular element
if it is self-adjoint and if its distribution is given by the
semicircle law: $$\tau
(x^k)=\frac{2}{\pi}\int_{-1}^{1}t^k\sqrt{1-t^2}dt\,\,\forall k\in
\mb{N}.$$ A family $(\mc{M}_i)_{i\in I}$ of unital $*$-subalgebras
of $\mc{M}$ is a free family if $\tau (x_k)=0$, $x_k\in
\mc{M}_{i_k}$ $\forall 1\leq k\leq p$, $i_1,\ldots ,i_p\in I$,
$i_1\not= i_2\not=\ldots\not= i_p$, $p\in \mb{N}$ imply $\tau
(x_1\ldots x_p)=0$. A family $(A_i)_{i\in I}$ of subsets
$A_i\subset \mc{M}$ is free if the family $(*$-alg$(\{1\}\cup
A_i))_{i\in I}$ is free. A free set $(s_i)_{1\leq i\leq m}\subset
\mc{M}$ consisting of semicircular elements is called a
semicircular system. If $(x_i)_{1\leq i\leq m}$ is free then $\chi
(x_1,\ldots ,x_m)=\chi (x_1)+\ldots +\chi (x_m)$ hence a finite
semicircular system has finite free entropy. The modified free
entropy dimension and the free entropy dimension of an $m$-tuple
of self-adjoint elements $(x_i)_{1\leq i\leq m}\subset\mc{M}$ are
$$\delta_0((x_i)_{1\leq i\leq m})=m+\limsup_{\omega\rightarrow
0}\frac{\chi ((x_i+\omega s_i)_{1\leq i\leq m}:(s_i)_{1\leq i\leq
m})}{|\log\omega |}$$ and respectively, $$\delta ((x_i)_{1\leq
i\leq m})=m+\limsup_{\omega\rightarrow 0}\frac{\chi ((x_i+\omega
s_i)_{1\leq i\leq m})}{|\log\omega |}\,,$$ where $(x_i)_{1\leq
i\leq m}$ and the semicircular system $(s_i)_{1\leq i\leq m}$ are
free. If $x_1\ldots ,x_m$ are free, then $$\delta_0((x_i)_{1\leq
i\leq m})=\delta ((x_i)_{1\leq i\leq m})=\sum_{i=1}^m\delta
(x_i).$$
Moreover, for a single self-adjoint element $x\in\mc{M}$ one has
$$\delta (x)=1-\sum_{s\in\mb{R}}(\mu (\{s\}))^2,$$ therefore $\delta
(x)=1$ if the distribution of $x$ has no atoms.
\section{Estimate of free entropy}\label{s2}
We obtain an estimate of the free entropy $\chi (x_1,\ldots ,x_m)$
for self-adjoint elements $x_1,\ldots ,x_m$ which can be
approximated in the $||\cdot ||_2$-norm by certain noncommutative
polynomials of degree $1$ in some of their variables. The proof of
Lemma \ref{l1} is based on the observation that in this case the
$c\times c$ matricial microstates of $x_1,\ldots ,x_m$ are
concentrated in some neighborhood of a linear subspace in
$\mc{M}_c^{sa}(\mb{C})$.
\begin{lema}\label{l1}
Let $x_1,\ldots , x_m$ be self-adjoint elements that generate a
$\mbox{\!I\!I}_1$-factor $(\mc{M},\tau )$. Assume that there exist
self-adjoint elements $m_j^{(l)},$ $z_k\in \mc{M}$ (for $1\leq
j\leq r+1$, $1\leq l\leq 2$, $1\leq k\leq 2v$), mutually
orthogonal projections $p_q\in \mc{M}$ (for $1\leq q\leq u$),
noncommutative polynomials
$\Phi_{ji}^{(l)}((p_q)_q,(z_k)_k)=\sum_{k=1}^{2v}\sum_{q,s=1}^u\mu_{q,s}^{(i,j,k,l)}
p_qz_kp_s$
(where $\mu_{q,s}^{(i,j,k,l)}$ are scalars), and $0<\omega
<\frac{1}{3}$ such that
\begin{displaymath}
\left|\left|x_i-\frac{1}{2}\sum_{l=1}^{2}\sum_{j=1}^{r+1}(m_j^{(l)}
\Phi_{ji}^{(l)}((p_q)_q,(z_k)_k)
+\Phi_{ji}^{(l)}((p_q)_q,(z_k)_k)^*
m_j^{(l)})\right|\right|_2<\omega
\end{displaymath}
for all $1\leq i\leq m$. Then
\begin{equation}\label{eq2}
\chi (x_1,\ldots,x_m)\leq C(m,r,v,K)+(m-2r-2v-3)\log\omega ,
\end{equation}
where $C(m,r,v,K)$ is a constant depending only on $m$, $r$, $v$,
and $K=1+\max_{i,j,l}\{|| \Phi_{ji}^{(l)}((p_q)_q,(z_k)_k)||_2,
||x_i||, ||m_j^{(l)}||\}$.
\end{lema}
\begin{proof} For $R,\frac{1}{\epsilon} >0$ sufficiently large and
integer $p\geq 1$ consider
$(A_1,\ldots,A_m,(M_j^{(l)})_{j,l},(P_q)_q,(Z_k)_k)$, an arbitrary
element of the set of matricial microstates $\Gamma_R
(x_1,\ldots,x_m,(m_j^{(l)})_{j,l},(p_q)_q,(z_k)_k;p,c,\epsilon)$.
One can assume (see \cite{15}) that $||A_i||, ||M_j^{(l)}||,
||P_q||\leq K$. If $p$ is large and $\epsilon
>0$ is small enough, then
\begin{eqnarray}&&
\bigg|\bigg|A_i-\frac{1}{2}\sum_{l=1}^{2}\sum_{j=1}^{r+1}(M_j^{(l)}
\Phi_{ji}^{(l)}((P_q)_q,(Z_k)_k)\\\nonumber &&\hspace{1 cm}
+\Phi_{ji}^{(l)}((P_q)_q,(Z_k)_k)^*
M_j^{(l)})\bigg|\bigg|_2<\omega
\end{eqnarray}
for all $1\leq i\leq m$ and
$||\Phi_{ji}^{(l)}((P_q)_q,(Z_k)_k)||_2< K$ for all $i, j, l$.
Lemma 4.3 in \cite{14} implies that for any $\delta
>0$ there exist $p',c'\in\mb{N}, \epsilon_1>0$ such that if $c\geq
c'$ and if $(P_1,\ldots
,P_u)\in\Gamma_R((p_q)_q;p',c,\epsilon_1)$, then there exist
mutually orthogonal projections $Q_1,\ldots ,Q_u\in
\mc{M}_c^{sa}(\mb{C})$ such that $\mbox{rank}(Q_q)=[\tau (p_q)c]$
and $||P_q-Q_q||_2<\delta\forall 1\leq q\leq u$. If $\delta >0$ is
sufficiently small one has then for all $c\geq c'$ and for all
$1\leq i\leq m$,
\begin{eqnarray}&&
\bigg|\bigg|A_i-\frac{1}{2}\sum_{l=1}^{2}\sum_{j=1}^{r+1}(M_j^{(l)}
\Phi_{ji}^{(l)}((Q_q)_q,(Z_k)_k)\\\nonumber
 &&\hspace{1 cm}+\Phi_{ji}^{(l)}((Q_q)_q, (Z_k)_k)^*
M_j^{(l)})\bigg|\bigg|_2<\omega
\end{eqnarray}
and $||\Phi_{ji}^{(l)}((Q_q)_q,(Z_k)_k)||_2<K$ for all $i, j, l$.
Let $S_1,\ldots ,S_u\in \mc{M}_c^{sa}(\mb{C})$ be mutually
orthogonal projections, fixed, with each projection $S_q$ of rank
$[\tau (p_q)c]$. There exists then $U\in \mc{U}_c(\mb{C})$
such that
 $Q_q=U^*S_qU$ for all $1\leq q\leq u$ and one obtains
\begin{eqnarray}&&
\bigg|\bigg|UA_iU^*-\frac{1}{2}\sum_{l=1}^{2}\sum_{j=1}^{r+1}(B_j^{(l)}
\Phi_{ji}^{(l)}((S_q)_q,(T_k)_k)\\\nonumber &&\hspace{1 cm}
+\Phi_{ji}^{(l)}((S_q)_q,(T_k)_k)^*
B_j^{(l)})\bigg|\bigg|_2<\omega
\end{eqnarray} for
all $1\leq i\leq m$, where we denoted $B_j^{(l)}=UM_j^{(l)}U^*$,
$T_k=UZ_kU^*$. Let $\{U_a\}_{a\in A(c)}$ be a minimal $\gamma$-net
in $\mc{ U}_c(\mb{C})$ with respect to the $||\cdot ||$-norm.
According to a result of S. J. Szarek (\cite{12}),
$|A(c)|\leq\left(\frac{C}{\gamma } \right)^{c^2}$ for some
universal constant $C$. Consider also a minimal $\theta$-net
$\{V_b\}_{b\in B(c, K)}$ in $\{B\in \mc{M}_c^{sa}(\mb{C}):
||B||\leq K\}$, with respect to the same norm. It is easily seen
that Szarek's result implies $|B(c,K)|\leq
\left(\frac{CK}{\theta}\right)^{c^2+c}$. Since
$||UA_iU^*-U_aA_iU_a^*||_2<2K\gamma$ for some $a\in A(c)$ and
$||B_j^{(l)}-V_{b(j,l)}||<\theta$ for some $b(j,l)\in B(c, K)$, we
have
\begin{eqnarray}\label{ineq}&&
\bigg|\bigg|U_aA_iU_a^*-\frac{1}{2}\sum_{l=1}^{2}\sum_{j=1}^{r+1}(V_{b(j,l)}
\Phi_{ji}^{(l)}((S_q)_q,(T_k)_k)
+\Phi_{ji}^{(l)}((S_q)_q,\\\nonumber && \hspace{.5 cm}(T_k)_k)^*
V_{b(j,l)})\bigg|\bigg|_2\leq ||UA_iU^*-U_aA_iU_a^*||_2
+\bigg|\bigg|UA_iU^*-\frac{1}{2}\\\nonumber &&\hspace{.5 cm}
\cdot\sum_{l=1}^{2}\sum_{j=1}^{r+1}(B_j^{(l)}
\Phi_{ji}^{(l)}((S_q)_q,(T_k)_k)+\Phi_{ji}^{(l)}((S_q)_q,(T_k)_k)^*
B_j^{(l)})\bigg|\bigg|_2
\\\nonumber
&&\hspace{.5 cm}
+\sum_{l=1}^{2}\sum_{j=1}^{r+1}||B_j^{(l)}-V_{b(j,l)}||\cdot
||\Phi_{ji}^{(l)}((S_q)_q,(T_k)_k)||_2\\\nonumber &&\hspace{.5
cm}< 2K\gamma +\omega +2(r+1)K\theta =3\omega\,\forall 1\leq i\leq
m.
\end{eqnarray}
Choose $\gamma =\frac{\omega}{2K}$, $\theta
=\frac{\omega}{2(r+1)K}$, and define the function
$F=(F_i((T_k)_k))_i:(\mc{M}_c^{sa}(\mb{C}))^{2v}\rightarrow
(\mc{M}_c^{sa}(\mb{C}))^m$ by
\begin{eqnarray}&&
F_i((T_k)_k)=\frac{1}{2}\sum_{l=1}^{2}\sum_{j=1}^{r+1}U_a^*(V_{b(j,l)}
\Phi_{ji}^{(l)}((S_q)_q,(T_k)_k)\\\nonumber &&\hspace{1 cm}
+\Phi_{ji}^{(l)}((S_q)_q,(T_k)_k)^* V_{b(j,l)})U_a\,\forall 1\leq
i\leq m.
\end{eqnarray}
It follows from (\ref{ineq}) that the distance in the euclidian
norm from the microstate $(A_1,\ldots ,A_m)$ to the image of $F$
is less than or equal to $3\omega\sqrt{mc}$. The polynomials
$\Phi_{ji}^{(l)}$ are linear in $(T_k)_k$, hence the image of $F$
is a linear subspace in $(\mc{M}_c^{sa}(\mb{C}))^m$, of dimension
$d_F\leq 2vc^2$. Denote by $L_F(\omega ,c)$ the intersection of
this subspace with the ball of euclidian radius $(3\omega
+K)\sqrt{mc}$ and by $B_F(\omega ,c)$ the cartesian product of
$L_F(\omega ,c)$ with the ball of (euclidian) radius
$3\omega\sqrt{mc}$ in the orthogonal complement of the image of
$F$. The set of matricial microstates
$\Gamma_R(x_1,\ldots,x_m:(m_j^{(l)})_{j,l},(p_q)_q,(z_k)_k;p,c,\epsilon
)$ is contained in $\bigcup_FB_F(\omega ,c)$, hence
\begin{eqnarray}&&
\mbox{vol}_{mc^2}(\Gamma_R(x_1,\ldots,x_m:(m_j^{(l)})_{j,l},(p_q)_q,(z_k)_k;p,c,\epsilon
)) \\\nonumber && \hspace{.5 cm}
\leq\sum_F\mbox{vol}_{mc^2}(B_F(\omega ,c))\\\nonumber
&&\hspace{.5 cm} =\sum_F\mbox{vol}_{d_F}((3\omega
+K)\sqrt{mc})\cdot
\mbox{vol}_{mc^2-d_F}(3\omega\sqrt{mc})\\\nonumber &&\hspace{.5
cm} =\sum_F\frac{(\pi mc)^{\frac{d_F}{2}}(3\omega
+K)^{d_F}}{\Gamma(1+\frac{d_F}{2})}\cdot \frac{(\pi
mc)^{\frac{mc^2-d_F}{2}}(3\omega
)^{mc^2-d_F}}{\Gamma(1+\frac{mc^2-d_F}{2})}\\\nonumber
&&\hspace{.5 cm}\leq \left(\frac{2CK}{\omega }\right)^{c^2}
\cdot\left[\left(\frac{2(r+1)CK^2}{\omega
}\right)^{c^2+c}\right]^{2(r+1)}\\\nonumber &&\hspace{.5 cm}
\cdot\frac{(\pi mc)^{\frac{mc^2}{2}}(2K)^{2vc^2}(3\omega
)^{(m-2v)c^2}2^{mc^2} }{\Gamma(1+\frac{mc^2}{2})}\,.
\end{eqnarray}
After taking the limit after
$c,p,\frac{1}{\epsilon}\rightarrow\infty$ in the resulting upper
bound for
$\chi_R(x_1,\ldots,x_m:(m_j^{(l)})_{j,l},(p_q)_q,(z_k)_k;p,c,\epsilon)$,
eliminating $R$ as in the definition of free entropy, and
recalling that $\{x_1,\ldots,x_m\}$ is a system of generators, one
obtains
\begin{eqnarray}\label{eq3}
\chi (x_1,\ldots,x_m)&=&\chi(x_1,\ldots,x_m:(m_j^{(l)})_{j,l},(p_q)_q,(z_k)_k)\\
\nonumber &\leq& C(m,r,v,K)+(m-2r-2v-3)\log\omega\,.
\end{eqnarray}
\end{proof}
\section{Infinite multiplicity}\label{s3}
Let $\mc{P}$ be a von Neumann algebra. If $\mc{Q}\subset\mc{P}$ is
a subalgebra, then the normalizer of $\mc{Q}$ in $\mc{P}$ is by
definition the set $N_\mc{P}(\mc{Q})=\{u\in
\mc{P}:uu^*=u^*u=1,u\mc{Q}u^*=\mc{Q}\}$.
\begin{defi}(\cite{17}, \cite{3})\label{d1}
A Cartan subalgebra of a von Neumann algebra $\mc{P}$ is a maximal
abelian $*$-subalgebra (MASA) $\mc{A}\subset \mc{P}$ such that:
i) $\mc{A}$ is the range of a normal conditional expectation; ii)
the normalizer $N_\mc{P}(\mc{A})$ of $\mc{A}$ in $\mc{P}$
generates $\mc{P}$.
\end{defi}
If $\mc{N}$ is a type $\mbox{\!I\!I}_1$-factor, then the
representation $\mc{N}\subset \mc{B}(L^2(\mc{ N},\tau))$ ($\tau$
denotes the unique normalized trace on $\mc{ N}$) is the standard
form of $\mc{N}$. Let $J:L^2(\mc{ N},\tau)\rightarrow
L^2(\mc{N},\tau)$ be the modular conjugacy operator. We recall the
following Theorem due to J. Feldman and C. C. Moore:
\begin{teo}\label{t1}(\cite{3}, \cite{9}) Let $\mc{N}$ be a type
$\mbox{\!I\!I}_1$-factor. If $\mc{A}$ is a Cartan subalgebra of
$\mc{N}$, then the algebra $(\mc{A}\cup J\mc{A}J)''$ is maximal
abelian in $\mc{B}(L^2(\mc{N},\tau))$.
\end{teo}
Being a MASA, the algebra $(\mc{A}\cup J\mc{A}J)''$ has a cyclic
vector $\xi\in L^2(\mc{N},\tau)$ i.e.,
$\overline{\mbox{sp}}^{||\cdot ||_2}(\mc{A}\cup
J\mc{A}J)''\xi=L^2(\mc{N},\tau)$. With the usual identification of
$_{J\mc{A}J}L^2(\mc{N},\tau)$ with $L^2(\mc{N},\tau)_\mc{A}$, this
means that $\overline{\mbox{sp}}^{||\cdot ||_2}\mc{A}\xi \mc{A}=
L^2(\mc{N},\tau)$ that is, $\mc{A}$ has finite multiplicity $1$ in
$\mc{N}$.
\begin{defi}(\cite{2})\label{d2}
An abelian subalgebra $\mc{A}$ of a type $\mbox{\!I\!I}_1$-factor
$\mc{N}$ has {\it finite multiplicity $\leq v<\infty$} if there
exist $v$ vectors $\xi_1,\ldots ,\xi_v\in L^2(\mc{N},\tau)$ such
that
$$\overline{\mbox{sp}}^{||\cdot ||_2}(\mc{A}\xi_1\mc{A}+\ldots +\mc{A}\xi_v\mc{
A}) =L^2(\mc{N},\tau)$$ or equivalently, if
$_\mc{A}L^2(\mc{N},\tau)_\mc{A}$ is generated as an
$\mc{A}$-$\mc{A}$-bimodule by $v$ vectors from $L^2(\mc{
N},\tau)$. If $_\mc{A}L^2(\mc{N},\tau)_\mc{A}$ is not a finitely
generated $\mc{A}$-$\mc{A}$-bimodule, we say that $\mc{ A}$ has
infinite multiplicity.
\end{defi}
The multiplicity of $\mc{A}$ in $\mc{N}$ does not increase after
compressing with a projection $p\in \mc{ A}$:
\begin{lema}(\cite{2})\label{l3}
If $\mc{A}\subset \mc{N}$ has finite multiplicity $\leq v$ and
$p\in \mc{A}$ is an arbitrary projection, then $\mc{
A}_p=p\mc{A}\subset p\mc{N}p=\mc{N}_p$ has also finite
multiplicity $\leq v$.
\end{lema}
\begin{teo}\label{t3}
Let $(\mc{M},\tau )$ be a $\mbox{\!I\!I}_1$-factor generated by
the self-adjoint elements $x_1,\ldots ,x_m$. If $\mc{ N}\subset
\mc{M}$ is a subfactor with the integer part of the Jones index
$[\mc{M}:\mc{N}]$ equal to $r$ and if $\mc{A}\subset\mc{N}$ is an
abelian subalgebra of multiplicity $\leq v$, then
$$\delta_0(x_1,\ldots ,x_m)\leq 2r+2v+3\,.$$
\end{teo}
\begin{proof} We can assume from the beginning that $m>2r+2v+3$ since
$\delta_0(x_1,\ldots ,x_m)\leq m$ is always true (\cite{15}).
There exists a Pimsner-Popa basis (\cite{8}) $m_1,\ldots
,m_{r+1}\in \mc{M}$ such that
$$x=\sum_{j=1}^{r+1}m_jE_\mc{N}(m_j^*x)\,\forall x\in \mc{M},$$
where $E_\mc{N}:\mc{M}\rightarrow \mc{N}$ is the conditional
expectation from $\mc{M}$ onto $\mc{N}$. Denote the embedding
$\mc{N}\subset L^2(\mc{N},\tau)$ by $x\mapsto \hat x$ and let
$J:L^2(\mc{N},\tau)\rightarrow L^2(\mc{ N},\tau)$ be the modular
conjugacy operator defined by $J\left(\hat x\right)=\widehat
{x^*}$. Let $\xi_1,\ldots,\xi_v\in L^2(\mc{N},\tau)$ such that
$$\mc{A}\xi_1\mc{A}+\ldots +\mc{A}\xi_v\mc{A}$$
is a dense subset of $L^2(\mc{N},\tau)$. Eventually replacing
$\xi_i$ by
$\frac{1}{2}(\xi_i+J\xi_i)+\frac{1}{2\sqrt{-1}}(\xi_i-J\xi_i)\sqrt{-1}$
and regrouping, we can assume that there exist
$\eta_1,\ldots,\eta_{2v}\in L^2(\mc{N},\tau)^{sa}:= \{\xi\in
L^2(\mc{N},\tau):J\xi =\xi \}$ such that $\mc{
A}\eta_1\mc{A}+\ldots +\mc{A}\eta_{2v}\mc{A}$ is dense in
$L^2(\mc{N},\tau)$. Let $x_1,\ldots ,x_m$ be self-adjoint elements
of $\mc{M}$. Every element $E_\mc{N}(m_j^*x_i)\in \mc{ N}$ can be
approximated arbitrarily well in the $||\cdot ||_2$-norm by
elements of the form
$$\sum_{k=1}^{2v}\sum_{p=1}^{t}a_{p,k}^{(i,j)}\eta_kb_{p,k}^{(i,j)}$$
for some $a_{p,k}^{(i,j)},b_{p,k}^{(i,j)}\in \mc{A}$. Since $\mc{A}$ is
abelian, there exist an integer $u$ and projections
$p_1,\ldots,p_u$ of sum $1$ such that every $a_{p,k}^{(i,j)}$ and
$b_{p,k}^{(i,j)}$ can be approximated sufficiently well in the uniform
norm, by linear combinations of these projections. Moreover,
$\widehat {\mc{N}^{sa}}$ is dense in $L^2(\mc{ N},\tau)^{sa}$ so
one can find $z_1,\ldots,z_{2v}$ self-adjoint elements of $\mc{N}$
and scalars $\mu_{q,s}^{(i,j,k)}\in\mb{C}$ such that
$$\Psi_{ji}((p_q)_q,(z_k)_k)=\sum_{k=1}^{2v}\sum_{q,s=1}^{u}\mu_{q,s}^{(i,j,k)}p_qz_kp_s$$
is sufficiently close to $E_\mc{N}(m_j^*x_i)$ in the
$||\cdot||_2$-norm, for all indices $i,j$. In particular, one can arrange
for the norms $|| \Psi_{ji}((p_q)_q,(z_k)_k)||_2$ to be all
uniformly bounded by a constant $D$ depending only on the norms
$||m_j^*x_i||$. Therefore, every element $x_i$ can be approximated
arbitrarily well in the $||\cdot ||_2$-norm, by elements of the
form
$$\sum_{j=1}^{r+1}m_j\Psi_{ji}((p_q)_q,(z_k)_k).$$
Denote $m_j^{(1)}=\frac{1}{2}(m_j+m_j^*)$ and
$m_j^{(2)}=\frac{1}{2\sqrt{-1}}(m_j-m_j^*)$. It follows that every
element $x_i$ can be approximated arbitrarily well in the $||\cdot
||_2$-norm, by elements of the form
$$\sum_{l=1}^2\sum_{j=1}^{r+1}m_j^{(l)}\Phi_{ji}^{(l)}((p_q)_q,(z_k)_k),$$
where
$\Phi_{ji}^{(1)}((p_q)_q,(z_k)_k)=\Psi_{ji}((p_q)_q,(z_k)_k)=-\sqrt{-1}\Phi_{ji}^{(2)}(
(p_q)_q,(z_k)_k)$. Since $x_i=x_i^*\forall 1\leq i\leq m$, given $\omega >0$, every
element $x_i$ can ultimately be approximated in
the $||\cdot ||_2$-norm as
\[
\left|\left|x_i-\frac{1}{2}\sum_{l=1}^{2}\sum_{j=1}^{r+1}(m_j^{(l)}
\Phi_{ji}^{(l)}((p_q)_q,(z_k)_k)
+\Phi_{ji}^{(l)}((p_q)_q,(z_k)_k)^*
m_j^{(l)})\right|\right|_2<\omega\,.
\]
If $s_1,\ldots ,s_m$ is a semicircular system free from
$x_1,\ldots ,x_m$ then (\cite{15})
\begin{eqnarray}&&\chi ((x_i+\omega s_i)_{1\leq
i\leq m}:(s_i)_{1\leq i\leq m}) =\chi \big((x_i+\omega s_i)_{1\leq
i\leq m}:\\\nonumber && \hspace{.5 cm}(s_i)_{1\leq i\leq
m},(m_j^{(l)})_{j,l},(p_q)_q,(z_k)_k\big)\leq \chi
\big((x_i+\omega s_i)_{1\leq i\leq m}:\\\nonumber &&\hspace{.5 cm}
(m_j^{(l)})_{j,l},(p_q)_q,(z_k)_k\big)
\end{eqnarray} since
$(m_j^{(l)})_{j,l},(p_q)_q,(z_k)_k\subset\{x_i+\omega
s_i,s_i\,|\,1\leq i\leq m\}''$. Note that
\begin{eqnarray}&&\bigg|\bigg|x_i+\omega
s_i-\frac{1}{2}\sum_{l=1}^{2}\sum_{j=1}^{r+1}(m_j^{(l)}
\Phi_{ji}^{(l)}((p_q)_q,(z_k)_k)\\\nonumber &&\hspace{1 cm}
+\Phi_{ji}^{(l)}((p_q)_q,(z_k)_k)^*
m_j^{(l)})\bigg|\bigg|_2<2\omega
\end{eqnarray} for all $1\leq i\leq m$, hence the estimate of free entropy from Lemma
\ref{l1} implies
$$\chi ((x_i+\omega s_i)_{1\leq i\leq m}:(s_i)_{1\leq i\leq m})\leq C(m,r,v,K)+(m-2r-2v-3)
\log 2\omega ,$$
therefore
\begin{eqnarray}&&\delta_0(x_1,\ldots ,x_m)=m+\limsup_{\omega\rightarrow
0}\frac{\chi ((x_i+\omega s_i)_{1\leq i\leq m}:(s_i)_{1\leq i\leq m})}{|\log\omega |}\\
\nonumber && \hspace{.5 cm}\leq m+\limsup_{\omega\rightarrow
0}\frac{C(m,r,v,K)+(m-2r-2v-3)\log 2\omega}{|\log\omega
|}\\\nonumber && \hspace{.5 cm}=m-(m-2r-2v-3)=2r+2v+3\,.
\end{eqnarray}
\end{proof}
\begin{cor}\label{c2}
The subfactors $\mc{N}$ of finite index in the interpolated free
group factors $\mc{L}(\mb{F}_t)$, $1<t\leq \infty$, do not contain
abelian subalgebras of finite multiplicity.
\end{cor}
\begin{proof} Consider first the case $1<t<\infty$ and suppose that $\mc{
N}$ has an abelian subalgebra $\mc{A}$ of finite multiplicity $\leq
v$. For every projection $p\in \mc{A}$, $p\mc{A}$ is an
abelian subalgebra of multiplicity $\leq v$ in $p\mc{N}p$ (Lemma
\ref{l3}). Moreover (\cite{7}), $\left[\mc{L}(\mb{
F}_t)_p:p\mc{N}p\right]=\left[\mc{L}(\mb{F}_t):\mc{
N}\right] <\infty$. Eventually replacing $\mc{A}$ by a MASA in
$\mc{N}$, which contains $\mc{A}$ (and thus, is of finite
multiplicity $\leq v$ in $\mc{N}$), we can assume that $\mc{A}$ is
a MASA in $\mc{N}$, hence has no minimal projections. Therefore,
there exists a projection $p\in \mc{A}$ such that
$m=1+\frac{t-1}{\tau (p)^2}$ is a conveniently large integer
(i.e., $m>2r+2v+3$). Theorem \ref{t3} implies that the (modified)
free entropy dimension of any finite system of generators of
$\mc{L}(\mb{F}_t)_p\simeq \mc{L}(\mb{F}_m)$ (compression
formula in \cite{1}, \cite{10}) is $\leq 2r+2v+3$, and this is in
contradiction with the fact that $\mc{L}(\mb{F}_m)$ is
generated by a semicircular system with $\delta_0=m$ (\cite{14},
\cite{15}).

Suppose now $t=\infty$ and let $x_1,x_2,\ldots$ be an infinite
semicircular system that generates $\mc{L}(\mb{F}_\infty)$.
Making use of inequality (\ref{eq3}), one
obtains:
$$\chi\big(x_1,\ldots,x_m:(m_j^{(l)})_{j,l},(p_q)_q,(z_k)_k\big)<\chi (x_1,\ldots,x_m),$$
for some suitable elements.
Let $E_n$, $n\geq 1$, be the conditional expectation from
 $\mc{L}(\mb{F}_\infty)$ onto $\{x_1,\ldots ,x_n\}''$. The convergence
in distribution (\cite{14}, \cite{15}) implies the existence
of a integer $n>m$ such that
\begin{eqnarray}&&\chi\bigg(x_1,\ldots,x_m:\left( E_n\left(m_j^{(l)}\right)\right)_{j,l},
(E_n\left(p_q\right))_q,(E_n\left(z_k\right))_k\bigg)\\\nonumber
&& \hspace{1 cm}<\chi (x_1,\ldots,x_m).
\end{eqnarray}
One obtains then a contradiction:
\begin{eqnarray}&&\chi(x_1,\ldots,x_n)
=\chi\bigg(x_1,\ldots,x_n:\left(E_n\left(m_j^{(l)}\right)\right)_{j,l},
(E_n\left(p_q\right))_q,\\\nonumber &&\hspace{.5 cm}
(E_n\left(z_k\right))_k\bigg)
\leq\chi\bigg(x_1,\ldots,x_m:\left(E_n\left(m_j^{(l)}\right)\right)_{j,l},
(E_n\left(p_q\right))_q,\\\nonumber &&\hspace{.5 cm}
(E_n\left(z_k\right))_k\bigg)+\chi (x_{m+1},\ldots,x_n)<\chi
(x_1,\ldots,x_m)\\\nonumber &&\hspace{.5 cm} +\chi
(x_{m+1},\ldots,x_n)=\chi(x_1,\ldots,x_n).
\end{eqnarray}
\end{proof}
\begin{cor}\label{c1}
The interpolated free group subfactors (of finite index) do not contain
Cartan subalgebras.
\end{cor}
\begin{proof} With the result of J. Feldman and C. C. Moore (Theorem
\ref{t1}), every Cartan subalgebra is in particular an abelian
subalgebra of multiplicity $1$, the statement follows immediately
from Corollary \ref{c2}.
\end{proof}


\begin{thebibliography}{13}
\bibitem[Di]{17} Dixmier, J.: {\it Sous anneaux ab\'{e}liens maximaux dans les facteurs de
type fini}. Ann. of Math. 59 (1954), 279-286
 \bibitem[Dy1]{1} Dykema, K.: {\it Interpolated free group factors}. Pac. J. Math. 163
 (1994), 123-135
 \bibitem[Dy2]{2} Dykema, K.: {\it Two applications of free entropy}. Math.
 Ann. 308 (1997), 547-558
 \bibitem[FeMo]{3} Feldman, J. and Moore, C. C.: {\it Ergodic equivalence
 relations
 cohomology and von Neumann algebras}. Trans. Amer. Math. Soc. 234 (1977), 289-361
 \bibitem[Ge]{4} Ge, L.: {\it Applications of free entropy to finite von Neumann
 algebras}. Amer. J. Math. 119 (1997), 467-485
 \bibitem[GePo]{5} Ge, L. and Popa, S.: {\it On some decomposition properties for
factors of type $\mbox{\!I\!I}_1$}. Duke Math. J. 94 (1998),
79-101
 \bibitem[Ha]{6} Haagerup, U.: {\it An Example of a Non Nuclear C$^*$-algebra which has the
 Metric Approximation Property}. Invent. Math. 50 (1979), 279-293
 \bibitem[Jo]{7} Jones, V. F. R.: {\it Index for Subfactors}. Invent. Math. 72
 (1983), 1-25
 \bibitem[PiPo]{8} Pimsner, M. and Popa, S.: {\it Entropy and index for subfactors}.
Ann. Scient. Ec. Norm. Sup. 19 (1986), 57-106
 \bibitem[Po]{9} Popa, S.: {\it Notes on Cartan subalgebras in type $\mbox{\!I\!I}_1$ factors}.
 Math. Scand. 57 (1985), 171-188
 \bibitem[R\u a]{10} R\u adulescu, F.: {\it Random matrices, amalgamated free
 products and subfactors of the von Neumann algebra of a free group, of
 noninteger index}. Invent. Math. 115 (1994), 347-389
 \bibitem[Sz]{12} Szarek, S. J.: {\it Nets of Grassmann manifolds and orthogonal
 group}. Proceedings of Research Workshop on Banach Space Theory (Bor-Luh-Lin,
 ed.), The University of Iowa, June 29-31 (1981), 169-185
 \bibitem[\c St]{11} \c Stefan, M. B.: {\it The primality of subfactors of finite index in the
interpolated free group factors}. Proc. of the AMS 126 (1998), 2299-2307
\bibitem[Ta]{16} Takesaki, M.: {\it Theory of Operator Algebras, $\mbox{\!I\!I}$,
$\mbox{\!I\!I\!I}$}.
Encyclopaedia of Mathematical Sciences 125, Springer-Verlag (2003)
\bibitem[Vo1]{13} Voiculescu, D.: {\it Circular and semicircular systems and free
 product factors}. Operator Algebras, Unitary Representations, Enveloping
 Algebras, and Invariant Theory, Progress in Mathematics, Volume 92,
 Birkh\"{a}user, Boston (1990), 45-60
 \bibitem[Vo2]{14} Voiculescu, D.: {\it The analogues of entropy and of
 Fisher's information measure in free probability theory, $\mbox{\!I\!I}$}. Invent. Math.
 118 (1994), 411-440
 \bibitem[Vo3]{15} Voiculescu, D.: {\it The analogues of entropy and of
 Fisher's information measure in free probability theory, $\mbox{\!I\!I\!I}$: the
 absence of Cartan subalgebras}. G.A.F.A. Vol. 6, No. 1 (1996), 172-199
 \end{thebibliography}
\end{document}